\documentclass[3p,12pt]{amsart}
\usepackage{amsmath,amssymb}
\usepackage{supertabular}
\usepackage{ifthen}
\usepackage{alltt}
\usepackage{color}
\usepackage{etex}
\usepackage[margin=1.25in]{geometry}
\usepackage[latin1]{inputenc}
\usepackage{amsmath}
\usepackage{amsthm}
\usepackage[alphabetic]{amsrefs}
\usepackage{amsfonts}
\usepackage[all,cmtip]{xy}
\xyoption{rotate}
\usepackage{amsfonts}
\usepackage{amssymb}
\usepackage{amscd}
\usepackage{bbm}
\usepackage{pictexwd,dcpic}
\usepackage[T1]{fontenc}
\usepackage{lmodern}
\usepackage{setspace}
\usepackage{mathpartir}
\usepackage{enumerate,multicol}
\usepackage[svgnames]{xcolor}
\usepackage[colorlinks=true, citecolor=Dark Green, linkcolor=Maroon, urlcolor=Purple]{hyperref}
\usepackage{stmaryrd}
\usepackage{tikz-cd}
\def\Lrg{\mathcal{L}_{\text{rg}}}

\def\L{\mathcal{L}}

\def\eqdef{=_{\text{df}}}

\def\eqdef{=_{\text{df}}}

\def\catcomp{\circ}

\newcommand{\card}[1]{\vert #1 \vert}

\newcommand{\U}{\mathcal{U}}

\newcommand{\setdefinition}[2]{
\lbrace #1 \:\vert\: #2 \rbrace
}

\newcommand{\powerset}[1]{\mathcal{P}(#1)}

\newcommand{\fourpartdef}[8]
{
	\left\{
		\begin{array}{lll|}
			#1 & \mbox{if } #2 \\
			#3 & \mbox{if } #4 \\
			#5 & \mbox{if } #6 \\
			#7 & \mbox{if } #8 \\
		\end{array}
	\right.
}

\newcommand{\twopartdef}[4]
{
	\left\{
		\begin{array}{ll}
			#1 & \mbox{if } #2 \\
			#3 & \mbox{if } #4 
		\end{array}
	\right.
}

\newcommand{\threepartdef}[6]
{
	\left\{
		\begin{array}{lll}
			#1 & \mbox{if } #2 \\
			#3 & \mbox{if } #4 \\
			#5 & \mbox{if } #6 
		\end{array}
	\right.
}

\newcommand{\ob}[1]{\text{ob}\:#1}
\newcommand{\morph}[1]{\text{mor}#1}

\newcommand{\commtrianglemapsto}[6]{
\xymatrix{
#1 \ar@{|->}[r]^{#2} \ar@{|->}[rd]_{#6} &#3 \ar@{|->}[d]^{#4}\\
          &#5  \\
}
}

\newcommand{\natcommsquaremapsto}[9]{
\xymatrix{
#1(#2) \ar[rrr]^{#5_{#2}} \ar[ddd]_{#1(#3)} & & &#6(#2) \ar[ddd]^{#6(#3)}\\
 &#7 \ar@{|->}[r] \ar@{|->}[d] &#8 \ar@{|->}[d]\\
 &#9 \ar@{|->}[r]          &\lbrace #9 \rbrace \\
#1(#4) \ar[rrr]_{#5_{#4}} & &         &#6(#4)  \\
}
}

\newcommand{\functordef}[8]{
\xymatrix{
#1 \ar[rr]^{#2} & &#3  \\
#4 \ar[d]_{#5} \ar@/_20pt/[dd]_{#7 \catcomp #5} & & #2(#4) \ar[d]_{#2(#5)} \ar@/^20pt/[dd]^{#2(#7 \catcomp #5)}\\
#6 \ar[d]_{#7} & & #2(#6) \ar[d]_{#2(#7)} \\
#8 & &#2(#8) \\
}
}

\newcommand{\functorarrowaction}[8]{
\xymatrix{
#1 \ar[rr]^{#2} & &#3  \\
#4 \ar[d]_{#5} & &#2(#4) \ar[d]_{#2(#5)} &#7 \ar@{|->}[d] \\
#6 & & #2(#6) &#8
}
}

\newcommand{\contrafunctorarrowaction}[8]{
\xymatrix{
#1 \ar[rr]^{#2} & &#3  \\
#4 \ar[d]_{#5} & &#2(#4) &#7 \\
#6 & & #2(#6) \ar[u]_{#2(#5)} &#8 \ar@{|->}[u] 
}
}

\newcommand{\universalarrow}[9]{
\xymatrix{
#1 \ar@{-->}[dd]_{#2}   &#4 \ar@{-->}[dd]_{#5} & &#7 \ar[ll]_{#8} \ar[lldd]^{#9} \\
 							   & \\
#3  						   &#6  \\
}
}

\newcommand{\rightrightarrowsxy}[4]{
\xymatrix{
#1 \ar@/^/[rr]^{#2} \ar@/_/[rr]_{#3} & &#4
}
}

\newcommand{\rightleftarrowsxy}[4]{
\xymatrix{
#1 \ar@/^/[rr]^{#2} & &#4 \ar@/^/[ll]^{#3}
}
}

\definecolor{xdxdff}{rgb}{0.49019607843137253,0.49019607843137253,1.}
\definecolor{uuuuuu}{rgb}{0.26666666666666666,0.26666666666666666,0.26666666666666666}
\definecolor{zzttqq}{rgb}{0.6,0.2,0.}
\definecolor{qqqqff}{rgb}{0.,0.,1.}
\definecolor{ffqqqq}{rgb}{1.,0.,0.}

\newcommand\frakfamily{\usefont{U}{yfrak}{m}{n}}
\DeclareTextFontCommand{\textfrak}{\frakfamily}

\def\eqgap{.2ex}
\def\overgap{.4ex}
\def\inferrulerule{.2pt}

\newlength\rulelength
\newlength\toplength
\newlength\bottomlength

\newcommand\myinferrule[2]{%
  \stackMath%
  \setlength\bottomlength{\widthof{$#1$}}%
  \setlength\toplength{\widthof{$#2$}}%
  \ifdim\toplength>\bottomlength%
    \setlength\rulelength{\the\toplength}%
  \else%
    \setlength\rulelength{\the\bottomlength}%
  \fi%
  \mathrel{%
    \stackunder[\overgap]{%
      \stackon[\overgap]{%
        \stackanchor[\eqgap]%
          {\rule{\the\rulelength}{\inferrulerule}}%
        {\rule{\the\rulelength}{\inferrulerule}}%
      }{#2}%
    }{#1}%
  }%
}

\theoremstyle{plain}\newtheorem{lemma}{Lemma}[section]
\theoremstyle{plain}
\theoremstyle{plain}\newtheorem{theorem}[lemma]{Theorem}
\theoremstyle{plain}\newtheorem{problem}[lemma]{Problem}
\theoremstyle{remark}\newtheorem{construction}[lemma]{Construction}
\theoremstyle{plain}
\theoremstyle{plain}
\theoremstyle{definition}\newtheorem{definition}[lemma]{Definition}
\theoremstyle{remark} \newtheorem{example}[lemma]{Example}
\theoremstyle{remark} 
\theoremstyle{remark} \newtheorem{remark}[lemma]{Remark}
\theoremstyle{plain}
\theoremstyle{definition}\newtheorem*{notation}{Notation}
\theoremstyle{definition}
\theoremstyle{definition}
\theoremstyle{definition}
\theoremstyle{definition}
\theoremstyle{definition}\newtheorem*{convention}{Convention}

\usepackage{newclude}

\author{Dimitris Tsementzis}
\author{Matthew Weaver}

\newif\ifcomments\commentstrue

\begin{document}

\date{\today}

%

\newcommand{\vars}[1]{\vert #1 \vert}
\newcommand{\dep}[2]{\text{proj}_{#1}(#2)}
\newcommand{\depe}[2]{\text{\emph{proj}}_{#1}(#2)}
\def\FOLDStt{TT_{Sig}}
\def\FOLDS{Sig}
\def\Con{Con}
\def\Sort{SN}
\def\Var{Var}
\def\freshfor{\: \# \:}
\def\notboundin{\notin}
\newcommand\alphaequiv[2]{#1 \equiv_\alpha #2}
\newcommand\alphaclass[1]{[#1]}
\newcommand\reif[1]{\llbracket #1 \rrbracket}
\newcommand\subcomp[2]{#1 \cdot #2} 
\newcommand\cosieve[2]{#1 \downarrow #2}
\newcommand\oo[2]{#1 <_o #2}
\newcommand\om[2]{#1 <_m #2} 
\newcommand\weaken[2]{\textbf{wk}(#1,#2)}
\newcommand\shiftcon[2]{\textbf{shiftcon}(#1,#2)}
\newcommand\shift[2]{\textbf{shift}(#1,#2)}

\newcommand{\ottdrule}[4][]{{\displaystyle\frac{\begin{array}{l}#2\end{array}}{#3}\quad\ottdrulename{#4}}}
\newcommand{\ottusedrule}[1]{\[#1\]}
\newcommand{\ottpremise}[1]{ #1 \\}
\newenvironment{ottdefnblock}[3][]{ \framebox{\mbox{#2}} \quad #3 \\[0pt]}{}
\newenvironment{ottfundefnblock}[3][]{ \framebox{\mbox{#2}} \quad #3 \\[0pt]\begin{displaymath}\begin{array}{l}}{\end{array}\end{displaymath}}
\newcommand{\ottfunclause}[2]{ #1 \equiv #2 \\}
\newcommand{\ottnt}[1]{\mathit{#1}}
\newcommand{\ottmv}[1]{\mathit{#1}}
\newcommand{\ottkw}[1]{\mathbf{#1}}
\newcommand{\ottsym}[1]{#1}
\newcommand{\ottcom}[1]{\text{#1}}
\newcommand{\ottdrulename}[1]{\textsc{#1}}
\newcommand{\ottcomplu}[5]{\overline{#1}^{\,#2\in #3 #4 #5}}
\newcommand{\ottcompu}[3]{\overline{#1}^{\,#2<#3}}
\newcommand{\ottcomp}[2]{\overline{#1}^{\,#2}}
\newcommand{\ottgrammartabular}[1]{\begin{supertabular}{llcllllll}#1\end{supertabular}}
\newcommand{\ottmetavartabular}[1]{\begin{supertabular}{ll}#1\end{supertabular}}
\newcommand{\ottrulehead}[3]{$#1$ & & $#2$ & & & \multicolumn{2}{l}{#3}}
\newcommand{\ottprodline}[6]{& & $#1$ & $#2$ & $#3 #4$ & $#5$ & $#6$}
\newcommand{\ottfirstprodline}[6]{\ottprodline{#1}{#2}{#3}{#4}{#5}{#6}}
\newcommand{\ottlongprodline}[2]{& & $#1$ & \multicolumn{4}{l}{$#2$}}
\newcommand{\ottfirstlongprodline}[2]{\ottlongprodline{#1}{#2}}
\newcommand{\ottbindspecprodline}[6]{\ottprodline{#1}{#2}{#3}{#4}{#5}{#6}}
\newcommand{\ottprodnewline}{\\}
\newcommand{\ottinterrule}{\\[5.0mm]}
\newcommand{\ottafterlastrule}{\\}
\renewcommand{\ottdrule}[4][]{ {\displaystyle\frac{\begin{array}{c}#2\end{array} }{#3}\quad\ottdrulename{#4} } }
  
\newcommand{\ottmetavars}{
\ottmetavartabular{
 $ \ottmv{variables} ,\, \ottmv{x} ,\, \ottmv{z} $ &  \\
 $ \ottmv{sorts} ,\, \ottmv{A} $ &  \\
 $ \ottmv{index} ,\, \ottmv{i} $ &  \\
}}

\newcommand{\ottP}{
\ottrulehead{\Psi, \Phi}{::=}{\ottcom{Signatures}}\ottprodnewline
\ottfirstprodline{|}{ \bullet }{}{}{}{\ottcom{empty signature}}\ottprodnewline
\ottprodline{|}{\Psi  \ottsym{,}  \ottmv{A}  \ottsym{:}  \Gamma}{}{}{}{\ottcom{signature extension}}}

\newcommand{\ottG}{
\ottrulehead{\Gamma  ,\ \Delta}{::=}{\ottcom{Contexts}}\ottprodnewline
\ottfirstprodline{|}{ \bullet }{}{}{}{\ottcom{empty context}}\ottprodnewline
\ottprodline{|}{\Gamma  \ottsym{,}  \ottmv{x}  \ottsym{:}  \ottnt{X}}{}{}{}{\ottcom{context extension}}}

\newcommand{\otts}{
\ottrulehead{\sigma, \tau, \delta}{::=}{\ottcom{Substitutions}}\ottprodnewline
\ottfirstprodline{|}{ \epsilon }{}{}{}{\ottcom{empty substitution}}\ottprodnewline
\ottprodline{|}{\sigma  \ottsym{,}  \ottnt{t}}{}{}{}{\ottcom{substitution extension}}}

\newcommand{\ottt}{
\ottrulehead{\ottnt{t}}{::=}{\ottcom{Variables}}\ottprodnewline
\ottfirstprodline{|}{\ottmv{x}}{}{}{}{x \in Var}}

\newcommand{\ottty}{
\ottrulehead{\ottnt{X}}{::=}{\ottcom{Sorts}}\ottprodnewline
\ottprodline{|}{\ottmv{A} \, \sigma}{}{}{}{A \in \Sort}}

\newcommand{\ottgrammar}{\ottgrammartabular{
\ottP\ottinterrule
\ottG\ottinterrule
\otts\ottinterrule
\ottt\ottinterrule
\ottty\ottafterlastrule
}}

\newcommand{\ottdrulesigXXempty}[1]{\ottdrule[#1]{%
}{
 \bullet   \vdash}{%
{\ottdrulename{sig\_empty}}{}%
}}

\newcommand{\ottdrulesigXXext}[1]{\ottdrule[#1]{%
\ottpremise{\Psi  \vdash \quad  \Psi ; \  \Gamma  \vdash  \quad \ottmv{A} \, \notboundin \, \Psi}%
}{
\Psi  \ottsym{,}  \ottmv{A}  \ottsym{:}  \Gamma  \vdash}{%
{\ottdrulename{sig\_ext}}{}%
}}

\newcommand{\ottdefnsigXXwf}[1]{\begin{ottdefnblock}[#1]{$\Psi  \vdash$}{\ottcom{Well-Formed Signatures}}
\ottusedrule{\ottdrulesigXXempty{}}
\bigskip
\ottusedrule{\ottdrulesigXXext{}}
\end{ottdefnblock}}

\newcommand{\ottdrulectxXXempty}[1]{\ottdrule[#1]{%
}{
 \Psi ; \   \bullet   \vdash }{%
{\ottdrulename{ctx\_empty}}{}%
}}

\newcommand{\ottdrulectxXXext}[1]{\ottdrule[#1]{%
\ottpremise{ \Psi ; \  \Gamma  \vdash  \quad  \Psi ; \  \Gamma  \vdash  \ottnt{X}    \quad \ottmv{x} \, \notboundin \, \Gamma}%
}{
 \Psi ; \  \Gamma  \ottsym{,}  \ottmv{x}  \ottsym{:}  \ottnt{X}  \vdash }{%
{\ottdrulename{ctx\_ext}}{}%
}}

\newcommand{\ottdefnctxXXwf}[1]{\begin{ottdefnblock}[#1]{$ \Psi ; \  \Gamma  \vdash $}{\ottcom{Well-Formed Contexts}}
\ottusedrule{\ottdrulectxXXempty{}}
\bigskip
\ottusedrule{\ottdrulectxXXext{}}
\end{ottdefnblock}}

\newcommand{\ottdrulesubXXempty}[1]{\ottdrule[#1]{%
}{
\Psi  \vdash   \epsilon   \ottsym{:}  \Gamma  \Rightarrow   \bullet }{%
{\ottdrulename{sub\_empty}}{}%
}}

\newcommand{\ottdrulesubXXext}[1]{\ottdrule[#1]{%
\ottpremise{ \Psi  \vdash  \sigma  \ottsym{:}  \Gamma  \Rightarrow  \Delta  \quad  \Psi ; \  \Gamma  \vdash  \ottmv{z}  :  \sigma(\ottnt{X}) }%
}{
\Psi  \vdash   ( \sigma  \ottsym{,}  \ottmv{z} )   \ottsym{:}  \Gamma  \Rightarrow  \Delta  \ottsym{,}  \ottmv{x}  \ottsym{:}  \ottnt{X}}{%
{\ottdrulename{sub\_ext}}{}%
}}

\newcommand{\ottdefnsubsts}[1]{\begin{ottdefnblock}[#1]{$\Psi  \vdash  \sigma  \ottsym{:}  \Gamma  \Rightarrow  \Delta$}{\ottcom{Substitutions}}
\ottusedrule{\ottdrulesubXXempty{}}
\bigskip
\ottusedrule{\ottdrulesubXXext{}}
\end{ottdefnblock}}

\newcommand{\ottdruletypeXXvar}[1]{\ottdrule[#1]{%
\ottpremise{\ottsym{(}  \ottmv{x}  \ottsym{:}  \ottnt{X}  \ottsym{)} \, \in \, \Gamma \quad \Psi;\: \Gamma \vdash}%
}{
 \Psi ; \  \Gamma  \vdash  \ottmv{x}  :  \ottnt{X} }{%
{\ottdrulename{type\_var}}{}%
}}

\newcommand{\ottdruletypeXXsort}[1]{\ottdrule[#1]{%
\ottpremise{\ottsym{(}  \ottmv{A}  \ottsym{:}  \Delta  \ottsym{)} \, \in \, \Psi \quad \Psi  \vdash  \sigma  \ottsym{:}  \Gamma  \Rightarrow  \Delta}%
}{
 \Psi ; \  \Gamma  \vdash  \ottmv{A} \, \sigma    }{%
{\ottdrulename{type\_sort}}{}%
}}

\newcommand{\ottdefntyping}[1]{\begin{ottdefnblock}[#1]{$ \Psi ; \  \Gamma  \vdash  \ottnt{t}  :  \ottnt{X} $}{\ottcom{Term-in-Sort}}
\ottusedrule{\ottdruletypeXXvar{}}
\end{ottdefnblock}}

\newcommand{\ottdefntypingu}[1]{\begin{ottdefnblock}[#1]{$ \Psi ; \  \Gamma  \vdash  \ottnt{X}  :  \U $}{\ottcom{Well-Formed Sort}}
\ottusedrule{\ottdruletypeXXsort{}}
\end{ottdefnblock}}


\newcommand{\ottdefnctxXXfresh}[1]{\begin{ottdefnblock}[#1]{$\ottmv{x} \, \notin \, \Gamma$}{}
\end{ottdefnblock}}


\newcommand{\ottdefnsigXXfresh}[1]{\begin{ottdefnblock}[#1]{$\ottmv{A} \, \notin \, \Psi$}{}
\end{ottdefnblock}}


\newcommand{\ottdefnctxXXelem}[1]{\begin{ottdefnblock}[#1]{$\ottsym{(}  \ottmv{x}  \ottsym{:}  \ottnt{X}  \ottsym{)} \, \in \, \Gamma$}{}
\end{ottdefnblock}}


\newcommand{\ottdefnsigXXelem}[1]{\begin{ottdefnblock}[#1]{$\ottsym{(}  \ottmv{A}  \ottsym{:}  \Gamma  \ottsym{)} \, \in \, \Psi$}{}
\end{ottdefnblock}}

\newcommand{\ottdrulectxXXcong}[1]{\ottdrule[#1]{%
\ottpremise{\ottnt{x}  \ \# \  \ottnt{Y} \quad \ottnt{y}  \ \# \  \ottnt{X}}%
}{
\Gamma  \ottsym{,}  \ottnt{x}  \ottsym{:}  \ottnt{X}  \ottsym{,}  \ottnt{y}  \ottsym{:}  \ottnt{Y}  \ottsym{,}  \Delta  \equiv  \Gamma  \ottsym{,}  \ottnt{y}  \ottsym{:}  \ottnt{Y}  \ottsym{,}  \ottnt{x}  \ottsym{:}  \ottnt{X}  \ottsym{,}  \Delta}{%
{\ottdrulename{ctx\_swap}}{}%
}}

\newcommand{\ottdrulectxXXrefl}[1]{\ottdrule[#1]{%
}{
\Gamma  \equiv  \Gamma}{%
{\ottdrulename{ctx\_refl}}{}%
}}

\newcommand{\ottdrulectxXXsym}[1]{\ottdrule[#1]{%
\ottpremise{\Gamma  \equiv  \Delta}%
}{
\Delta  \equiv  \Gamma}{%
{\ottdrulename{ctx\_sym}}{}%
}}

\newcommand{\ottdrulectxXXtrans}[1]{\ottdrule[#1]{%
\ottpremise{\Gamma  \equiv  \Delta \quad \Delta  \equiv  \Lambda}%
}{
\Gamma  \equiv  \Lambda}{%
{\ottdrulename{ctx\_trans}}{}%
}}

\newcommand{\ottdefnctxXXeq}[1]{\begin{ottdefnblock}[#1]{$\Gamma  \equiv  \Delta$}{\ottcom{Congruence on Contexts}}
\ottusedrule{\ottdrulectxXXcong{}}
\bigskip
\ottusedrule{\ottdrulectxXXrefl{}}
\bigskip
\ottusedrule{\ottdrulectxXXsym{}}
\bigskip
\ottusedrule{\ottdrulectxXXtrans{}}
\end{ottdefnblock}}

\newcommand{\ottdrulesigXXcong}[1]{\ottdrule[#1]{%
\ottpremise{\ottmv{A}  \ \# \  \Delta \quad \ottmv{B}  \ \# \  \Gamma}%
}{
\Psi  \ottsym{,}  \ottmv{A}  \ottsym{:}  \Gamma  \ottsym{,}  \ottmv{B}  \ottsym{:}  \Delta  \ottsym{,}  \Phi  \equiv  \Psi  \ottsym{,}  \ottmv{B}  \ottsym{:}  \Delta  \ottsym{,}  \ottmv{A}  \ottsym{:}  \Gamma  \ottsym{,}  \Phi}{%
{\ottdrulename{sig\_swap}}{}%
}}

\newcommand{\ottdrulesigXXrefl}[1]{\ottdrule[#1]{%
}{
\Psi  \equiv  \Psi}{%
{\ottdrulename{sig\_refl}}{}%
}}

\newcommand{\ottdrulesigXXsym}[1]{\ottdrule[#1]{%
\ottpremise{\Psi  \equiv  \Phi}%
}{
\Phi  \equiv  \Psi}{%
{\ottdrulename{sig\_sym}}{}%
}}

\newcommand{\ottdrulesigXXtrans}[1]{\ottdrule[#1]{%
\ottpremise{\Psi  \equiv  \Phi \quad \Phi  \equiv  \Xi}%
}{
\Psi  \equiv  \Xi}{%
{\ottdrulename{sig\_trans}}{}%
}}

\newcommand{\ottdefnsigXXeq}[1]{\begin{ottdefnblock}[#1]{$\Psi  \equiv  \Phi$}{\ottcom{Congruence on Signatures}}
\ottusedrule{\ottdrulesigXXcong{}}
\bigskip
\ottusedrule{\ottdrulesigXXrefl{}}
\bigskip
\ottusedrule{\ottdrulesigXXsym{}}
\bigskip
\ottusedrule{\ottdrulesigXXtrans{}}
\end{ottdefnblock}}

\newcommand{\ottdefnsRules}{
\ottdefnsigXXwf{}\ottdefnctxXXwf{}\ottdefnsubsts{}\ottdefntyping{}\ottdefntypingu{}
\ottdefnctxXXeq{}
\ottdefnsigXXeq{}
}

\newcommand{\ottdefnss}{
\ottdefnsRules
}

\newcommand{\ottall}{\ottmetavars\\[0pt]
\ottgrammar\\[5.0mm]
\ottdefnss}

\title{Finite Inverse Categories as Signatures}

\begin{abstract}
We define a simple dependent type theory and prove that its well-formed types correspond exactly to finite inverse categories.
\end{abstract}

\maketitle

\section{Introduction}

A finite inverse category (\emph{fic}) is a finite, skeletal category with no non-identity endomorphisms.
Fics are interesting from the point of view of dependent type theory because they correspond to collections of dependently-typed data.
The basic idea is well-illustrated by the following simple example.

\begin{example}\label{runningexample}
Consider the following fic

\hspace{5cm} \xymatrix{
&I \ar[d]_{i} \\
\Lrg \eqdef &A \ar@/^/[d]^{d} \ar@/_/[d]_{c} \\
&O
}

\noindent subject to the relation $di=ci$. $\Lrg$ can be thought of as a signature useful for formalizing reflexive graphs. The data in $\Lrg$ corresponds to the following collection of data in dependent type theory, writing $\U$ for a universe of types:
\begin{align*}
O &\colon \U \\
A &\colon O \times O \to \U \\ 
I &\colon (\Sigma \:(x \colon O) \:A(x,x)) \to \U
\end{align*}
or, equivalently, given the appropriate type constructors, to the type
\[
\Sigma \:(O \colon \U)\:(A \colon O \times O \to \U)\: ((\Sigma \:(x \colon O) \:A(x,x)) \to \U)
\]
\end{example}

We can therefore understand fics as a syntax that captures the ``nested $\Sigma$-types'' that one usually refers to as a dependently-typed signature.
However, if one wishes to use fics as a syntax in this manner, then one should be able to come up with an inductive process that produces any such fic, similar to how one can describe all first-order signatures as (recursively enumerable) sets of symbols with possibly associated arities, or to how one defines the well-formed terms of a type theory. 

In this paper, we carry out this task.
To that end we define a certain simple dependent type theory, which manages to produce exactly the dependently-typed signatures that correspond to fics. Another way to look at this formal system is that it encodes an (quotient-)inductive-inductive definition of the data type of all fics, similar to what is done in \cite{altenkirch16}. 

Both perspectives are important: the former allows us to work ``externally'' in order to define interpretations of fics in dependent type theory (as e.g. described in \cite{TsemHMT}), whereas the latter allows us to work ``internally'' and implement such an interpretation inside (extended versions of) type theory, as we intend to do in \cite{TsemWeavTTI}.

However, this paper is independent of either of these perspectives. Its aim is to simply define the relevant dependent type theory and prove (``externally'') that (a certain class of) well-formed data of this type theory correspond precisely to finite inverse categories as usually defined.

\subsection*{Related Work} The observation that finite inverse categories can play the role of a dependently-typed syntax was perhaps first made by Makkai in his work on First Order Logic With Dependent Sorts (FOLDS) \cite{MFOLDS}, with Cartmell's work \cite{Cart86} as an important precursor.
The connection between FOLDS and MLTT/HoTT has been pursued in several places, e.g. in \cite{TsemHMT, TsemHSIP, uniFOLDS} and pursuing connections of this sort are the primary motivation for our work here. 
In particular, we should single out Palmgren's work in \cite{PalmgrenFOLDS} where he develops a general notion of dependently-typed first-order logic, independent from our own, connecting some special cases to FOLDS, and thereby to fics. 
In particular, to the extent that his Theorems 4.1 and 4.3 overlap with our own Theorem \ref{sigbijfic}, he should certainly be given priority, although we arrived at our results independently (and for seemingly independent reasons).
Finite inverse categories can also be understood as special cases of Reedy categories, and so in that form they have been related to homotopy (type) theory, as explained for example in \cite{shulman2015reedy}. 

\subsection*{Outline} In Section \ref{tt}, we define the syntax and rules of the type theory $\FOLDStt$. In Section \ref{ttfic} we prove that the well-formed signatures of $\FOLDStt$ correspond exactly to finite inverse categories as usually defined (Theorem \ref{sigbijfic}).

\subsection*{Acknowledgments}
The first-named author was partially supported by NSF DMS-1554092 (P.I. Harry Crane).

\section{The Theory $\FOLDStt$}\label{tt}

We work in an extensional set theory, e.g. ZF.
We fix (countably infinite, disjoint) sets $\Var$ of \emph{variables} and $\Sort$ of \emph{sort names}.
We reserve the symbol $=$ for equality between elements of these sets, their combinations and equivalence classes. 
We will also assume that $\Var$ and $\Sort$ are equipped with some mechanism that allows us to produce, given any finite subset, a \emph{fresh} variable not in that subset. This can be achieved e.g. by regarding $\Var$ and $\Sort$ as nominal sets in the sense of \cite{PittsNS}.

\begin{definition}[Syntax of $\FOLDStt$]
The grammar of $\FOLDStt$ consists of \emph{signatures}, \emph{contexts}, \emph{substitutions} and \emph{variables} and \emph{sorts}, defined as follows: 

\ottgrammar
\end{definition}

\begin{remark}
If one wants to take the point of view that what we are defining here is a (non-standard) kind of dependent type theory, then our ``variables'' are the terms and our ``sorts'' are the types of this type theory.
\end{remark}

\begin{notation}
We will generally denote a substitution $(\dots((\sigma, x), y), \dots, z)$ as $(\sigma,x,y,z)$ and when all variables are explicit simply as a list of variables $(x,y,\dots,z)$.
\end{notation}

\begin{definition}[Judgments of $\FOLDStt$]
We have the following forms of judgment built out of our grammar, displayed together with their intended meaning:

\begin{tabular}{l | r}
$\Psi \vdash$ & ``$\Psi$ is a well-formed signature'' \\
 $\Psi; \Gamma \vdash$ &  ``$\Gamma$ is a well-formed context in signature $\Psi$'' \\
$\Psi; \Gamma \vdash X$ & ``$X$ is a sort in context $\Gamma$ and signature $\Psi$'' \\
$\Psi; \Gamma \vdash t \colon X$ & ``$t$ is a term of type $X$ in context $\Gamma$ and signature $\Psi$'' \\
$\Psi \vdash \sigma \colon \Gamma \Rightarrow \Delta$ & ``$\sigma$ is a context morphism from $\Gamma$ to $\Delta$ in signature $\Psi$'' \\
$\Psi \equiv \Phi$ & ``$\Psi$ and $\Phi$ are equivalent signatures'' \\
$\Gamma \equiv \Delta$ & ``$\Gamma$ and $\Delta$ are equivalent contexts''
\end{tabular}

\end{definition}

\begin{remark}
The judgment $\Psi; \Gamma \vdash X$ is our only ``typing'' judgment. Its intended meaning is that the sort $X$ is a well-formed type (in some universe $\U$).
\end{remark}




\begin{notation}
We use the following notation for stating fact about our (raw) syntax:
\begin{itemize}
\item $A \notboundin \Psi$ iff $A$ has not previously been bound in $\Psi$
\item $x \notboundin \Gamma$ iff $x$ has not previously been bound in $\Gamma$
\item $A \freshfor \Gamma$ iff $A$ is fresh for the context $\Gamma$
\item $x \freshfor X$ iff $x$ is fresh for the sort $X$
\item $(A \colon \Delta) \in \Psi$ if $A \colon \Delta$ is an element of the list $\Psi$
\item $(x \colon X) \in \Gamma$ if $x \colon X$ is an element of the list $\Gamma$
\end{itemize}
\end{notation}

\begin{remark}
The distinction between $\freshfor$ and $\notin$ is essential, since $\freshfor$ invariant with respect to $\alpha$-equivalence of expressions, whereas $\notin$ is not.
\end{remark}

\begin{definition}[Rules of $\FOLDStt$]
The rules of $\FOLDStt$ are as follows, grouped according the form of the judgment in their conclusion:

\bigskip

\ottdefnss

\end{definition}

In order for the symbol $\sigma (X)$ in rule \ottdrulename{sub\_ext} to make sense we require the following auxiliary definitions, carried out simultaneously with the derivation of well-formed contexts and substitutions, and mutually one to the other.

\begin{definition}[Variables in a context]\label{vars}
We define a function $\vars{-}$
 simultaneously with the definition of well-formed contexts in $\FOLDStt$ as follows:
\begin{itemize}
\item $\vars{\bullet} \eqdef \varnothing$
\item $\vars{\Gamma, x \colon X} \eqdef \vars{\Gamma} \amalg \lbrace x \rbrace$
\end{itemize} 
\end{definition}

\begin{definition}[Substitution as a function]
We define for every well-formed substitution $\sigma \colon \Gamma \Rightarrow \Delta$, a function $\reif{\sigma} \colon \vars{\Delta} \rightarrow \vars{\Gamma}$ by induction on $\sigma$ as follows:
\begin{itemize}
\item $\reif{\epsilon}$ is the unique arrow $\varnothing \rightarrow \vars{\Gamma}$
\item Assuming $\reif{\sigma}$ has been defined for $\sigma \colon \Gamma \Rightarrow \Delta$ we define:
\begin{align*}
\reif{\sigma,z} \colon \vars{\Delta, x \colon X} &\rightarrow \vars{\Gamma} \\
\reif{\sigma,z}(y) &= \twopartdef{\reif{\sigma}(y)}{y \in \vars{\Delta}}{z}{ y = x}
\end{align*}
\end{itemize}
\end{definition}

\begin{definition}[Composition]\label{subcomposition}
Given well-formed substitutions $\sigma \colon \Gamma \Rightarrow \Delta$ and $\tau \colon \Delta \Rightarrow \Sigma$ we define their \textbf{composite} $\subcomp{\sigma}{\tau}$ (diagrammatic ordering!) by case analysis as follows:
\begin{itemize}
\item $\subcomp{\sigma}{\epsilon} \eqdef \epsilon$
\item $\subcomp{\sigma}{(\tau,z)} \eqdef (\subcomp{\sigma}{\tau}), \reif{\sigma} (z)$
\end{itemize}
\end{definition}

With these definition we define the action of substitutions on well-formed sorts, as it appears in the rule \ottdrulename{sub\_ext}.

\begin{definition}
$\sigma (A\:\tau) \eqdef A\:\subcomp{\sigma}{\tau}$ 
\end{definition}

The following lemmas will be useful to us.

\begin{lemma}\label{reifassoc}
$\reif{\subcomp{\sigma}{\tau}} = \reif{\tau} \catcomp \reif{\sigma}$
\end{lemma}
\begin{proof}
Immediate by induction on the structure of $\tau$.
\end{proof}

\begin{lemma}\label{reifcomp}
For any derivable $\Psi \vdash \sigma \colon \Gamma \Rightarrow \Delta$, $\Psi \vdash \tau \colon \Delta \Rightarrow \Sigma$, if $(A \colon \Sigma) \in \Psi$ and $(y \colon A\:\tau) \in \Delta$ then we have $(\reif{\tau}(y) \colon A\:\subcomp{\sigma}{\tau}) \in \Gamma$.
\end{lemma}
\begin{proof}
Since $(y \colon A\:\sigma) \in \Delta$, then $\Gamma \vdash \reif{\tau}(y) \colon (A\:\subcomp{\sigma}{\tau})$ which, on the assumption that everything is well-formed, is only the case if $(\reif{\tau}(y) \colon A\:\subcomp{\sigma}{\tau}) \in \Gamma$.
\end{proof}

We want to regard the derivable signatures and contexts of $\FOLDStt$ only up to $\alpha$-equivalence, which means that we want to regard two signatures (resp. contexts) that differ only in the consistent renaming of the variables that have been introduced in their derivation tree at every binding instance of a variable.
This can be achieved in a number of ways, e.g. by regarding $\Var$ and $\Sort$ as nominal sets, and we will assume that we have one such way of expressing such a notion of $\alpha$-equivalence.

\begin{definition}
Two derivable signatures $\Psi, \Phi$ (resp. contexts $\Gamma, \Delta$) are \textbf{congruent} if $\Psi \equiv \Phi$ (resp. $\Gamma \equiv \Delta$) is derivable. We say that $\Psi$ is \textbf{isomorphic} to $\Phi$ if $\Psi$ and $\Phi$ are either congruent or $\alpha$-equivalent, and similarly for contexts. We write $[\Phi]$ or $[\Gamma]$ for the isomorphism class of a signature $\Phi$ or a context $\Gamma$.

We write $\Con$ for the set of isomorphism classes of well-formed contexts of $\FOLDStt$:
\[
Con \eqdef \setdefinition{\alphaclass{\Gamma}}{\Psi; \Gamma \vdash \text {is derivable in $\FOLDStt$ for some $\Psi$}}
\]

We write $\FOLDS$ for the set of isomorphism classes of well-formed signatures of $\FOLDStt$:
\[
\FOLDS \eqdef \setdefinition{\alphaclass{\Psi}}{\Psi \vdash \text {is derivable in $\FOLDStt$}}
\]
\end{definition}

\begin{remark}
We will refer to elements of $\FOLDS$ and $Con$ through representatives of the isomorphism classes that comprise them, as is customary.
\end{remark}


\section{$\FOLDStt$ signatures as finite inverse categories}\label{ttfic}

\begin{definition}
A \textbf{finite inverse category} is a category $\L$ such that
\begin{enumerate}
\item $\card{\morph{\L}}$ is finite.
\item $\L$ has no non-identity endomorphisms.
\item $\L$ is skeletal.
\end{enumerate} 
We write $FIC$ for the set of all isomorphism classes of finite inverse categories (i.e. the objects of a chosen skeleton of the category of finite inverse categories).
We denote by $\varnothing$ the empty finite inverse category.
\end{definition}

We now want to prove that the elements of $\FOLDS$ correspond exactly to finite inverse categories. 
We illustrate the basic idea behind the correspondence between $\FOLDS$ and $FIC$ with an example.

\begin{example}\label{illustration}
Let $\Lrg$ be the fic from Example \ref{runningexample}. Its corresponding signature (in the sense of $\FOLDStt$) is given by
\[
\Psi_{rg} \eqdef O \colon \bullet, A \colon (c \colon O, d \colon O), I \colon (x \colon O, i \colon A(x,x))
\]
Intuitively, the names $O,A,I$ correspond to the objects of $\Lrg$, the variables $c,d,i$ the arrows of $\Lrg$ and the substitution $(x,x) \colon (x \colon O) \Rightarrow (c \colon O, d \colon O)$ encodes the relation $di=ci$ in $\Lrg$ by thinking of $x$ as representing the composites $di$ and $ci$ (which are equal).
We can then show that $\Psi_{rg}$ is derivable, i.e. that $\Psi_{rg}$ is indeed an element of $\FOLDS$.
Similarly, we want to argue that any $\Psi \in \FOLDS$ gives rise to a fic whose objects are given by the sort names in $\Psi$, arrows by the variables in $\Psi$ and relations (if any) by the substitutions attached to sorts in $\Psi$.
\end{example}


\begin{notation}
For $\L \in FIC$ and $K \in \ob{\L}$ and $f \colon K \rightarrow K'$ we write $K_f$ for $K'$.
We write $\cosieve{K}{\L}$ for the total cosieve on $K$ in $\L$, i.e. the set of all non-identity maps with $K$ as their domain.


\end{notation}

%
%
%
%

We now define certain useful functions that allow us to extract out of the derivation trees for well-formed signatures certain information about the variables in them.

\begin{definition}
We define a function $\vars{-}_{(-)} \colon \Con \times \Sort \rightarrow \powerset{\Var}$ simultaneously with the definition of well-formed contexts in $\FOLDStt$ as follows:
\begin{align*}
\vars{\bullet}_X &= \varnothing \\
\vars{\Gamma, y \colon Y}_X &= \twopartdef{\vars{\Gamma}_X}{X \neq Y}{\vars{\Gamma}_X\amalg \lbrace y \rbrace}{X=Y}
\end{align*}
\end{definition}

\begin{definition}
For every derivable expression $\Psi; \Gamma \vdash$ and for all $(A \colon \Delta) \in \Psi$ we define a partial function
\[
\text{proj}^{\Psi, \Gamma, A, \Delta}_{(-)}{(-)} \colon \vars{\Delta} \times \vars{\Gamma} \rightarrow \vars{\Gamma}
\] 
simultaneously with the definition of well-formed signatures and with the definition of $\vars{-}$, as follows:
\[
\text{proj}^{\Psi, \Gamma, A, \Delta}_x(y) = \twopartdef{\reif{\sigma}(x)}{(y\colon A \sigma) \in \Gamma \text{ for } \sigma \colon \Gamma \Rightarrow \Delta}{\mathtt{fail}}{\text{otherwise}}
\]
\end{definition}

\begin{notation}
We will omit superscripts and write simply $\dep{x}{y}$.
\end{notation}

\begin{example}
Consider the signature $\Psi_{rg}$ from Example \ref{illustration} and using the same notation as there let $\Delta \eqdef c \colon O, d \colon O$ and $\Gamma \eqdef x \colon O, i \colon A (x,x)$. Then we have $\vars{\Delta} = \lbrace c,d \rbrace$ and $\vars{\Gamma} = \lbrace x ,i \rbrace$ and the substitution $(x,x)$ is given by the function $\reif{(x,x)}$ which sends both $c$ and $d$ to $x$. Then we have
$
\dep{c}{i} = \reif{(x,x)}(c) = x
$
and 
$
\dep{d}{i} = \reif{(x,x)}(d) = x
$.
This means that the ``dependency of $i$ at position $c$ (resp. $d$) is $x$'' which is exactly the intended meaning.
\end{example}

\begin{lemma}\label{depassociative}
Let $\Psi;\Delta$, $\Psi; \Gamma$, $\Psi ; \Sigma$ be derivable, $\sigma \colon \Gamma \Rightarrow \Delta$, $\tau \colon \Sigma \Rightarrow \Gamma$ be substitutions, $(A \colon \Delta) \in \Psi$ and $(y \colon A \, \sigma) \in \Gamma$, $(z \colon A \: \subcomp{\sigma}{\tau}) \in \Sigma$. 
Then $\depe{\depe{x}{y}}{z}$ and $\depe{x}{\depe{y}{z}}$ are defined and $$\depe{\depe{x}{y}}{z}=\depe{x}{\depe{y}{z}}$$
\end{lemma}
\begin{proof}
That $\dep{\dep{x}{y}}{z}$ and $\dep{x}{\dep{y}{z}}$ are defined is immediate from the assumption of the derivability of the given judgments. Hence, we have:
\begin{align*}
\dep{x}{\dep{y}{z}} &= \dep{x}{\reif{\tau}(y)} \tag{$z \colon A\:\subcomp{\sigma}{\tau}$} \\
			     &= \reif{\subcomp{\tau}{\sigma}} (x) \tag{Lemma \ref{reifcomp}} \\
			     &= \reif{\tau} (\reif{\sigma} (x)) \tag{Lemma \ref{reifassoc}} \\
			     &=\dep{\reif{\sigma}(x)}{z} \tag{$z \colon A\:\subcomp{\sigma}{\tau}$} \\
			     &=\dep{\dep{x}{y}}{z} \tag{$y \colon A\:\sigma$}	
\end{align*}
\end{proof}

\begin{remark}
Lemma \ref{depassociative} expresses the ``associativity'' of $\dep{(-)}{-}$ and we will use this property crucially in order to define an associative composition operation on the fic we will extract from a well-formed signature.
\end{remark}

We are now ready to construct functions between $FIC$ and $\FOLDS$ and prove that they are inverse to each other. It is helpful for these constructions to fix the following convention, which allows us to freely switch between symbols for objects and arrows of a fic and symbols for sorts and variables in $\FOLDStt$.

\begin{convention}
We let $FIC$ be the set of isomorphism classes of finite inverse categories $\L$ such that $\ob{\L} \subset \Sort$ and the non-identity morphisms in $\L$ are a subset of $\Var$. Furthermore, we assume fixed a formal symbol $1$ such that all the identity morphisms for an object $K$ in a fic are given by $1_K$.
\end{convention}


\begin{definition}\label{orderings}
For any $\L \in FIC$ we will assume fixed a partial order $\oo{}{}$ on $\ob{\L}$ such that $\exists f \colon K' \rightarrow K \Rightarrow \oo{K}{K'}$, and for any $K \in \ob{\L}$ a partial order $\om{}{}$ in $\cosieve{K}{\L}$ such that $h=g \catcomp f \Rightarrow \om{f}{h}$.
\end{definition}

\begin{remark}
Since both sets on which orders are defined in Definition \ref{orderings} are finite, such orders will always exist.
\end{remark}

\begin{definition}
Let $\L \in FIC$ and $K \in \ob{\L}$. Let $f_1,\dots, f_n$ be the ordered morphisms in $\cosieve{K}{\L}$ and for each $i=1,\dots,n$ let $p_{i1},\dots,p_{im_i}$ be the ordered morphisms in $\cosieve{K_{f_i}}{\L}$. We then define 
\[
T_K \eqdef (f_1 \colon K_{f_1} (p_{11}f_1, \dots, p_{1m_1}f_1),\dots,f_n \colon K_{f_n} (p_{n1}f_n, \dots, p_{nm_{n}}f_n))
\]
where the ordering is of the kind described in Definition \ref{orderings}, and the $(p_{i1}f_i, \dots, p_{im_i}f_i)$ are (lists of variables respresenting) substitutions.
\end{definition}


\begin{problem}
To construct a function $R \colon \FOLDS \rightarrow FIC$.
\end{problem}
\begin{construction}\label{fictosig}
We proceed by induction on the derivations of signatures $\Psi$.
If $\Psi \equiv \bullet$ then we define $R(\bullet) \eqdef \varnothing$.
Now assume $R(\Psi)$ is defined and let $\Phi \equiv \Psi, K \colon \Gamma$.
We will define a finite inverse category $R(\Phi)$ as follows:
\begin{itemize}
\item $\ob{R(\Phi)} = \ob{R(\Psi)} \amalg \lbrace K \rbrace$
\item 
\[
\morph{R(\Phi)}(X,Y) = \fourpartdef{\morph{R(\Psi)}(X,Y)}{ X,Y \neq K}{1_K}{X=Y=K}{\varnothing}{Y=K,X\neq K}{\vars{\Gamma}_Y}{X=K,Y\neq K}
\]
\item Let $f \in R(\Phi)(X,Y)$ and $g \in R(\Phi)(Y,Z)$. We define the composition $g \catcomp f \in R(\Phi)(X,Z)$ as follows:
\[
g \catcomp f = \threepartdef{g \catcomp_{R(\Psi)}f}{X,Y, Z \neq K}{g}{f=1_K}{\dep{g}{f}}{X=K, Y,Z \neq K}
\]
where $\catcomp_{R(\Psi)}$ is the composition operation of the category $R(\Psi)$ and $g$ in $\dep{g}{f}$ is identified with the variable from which it has been defined (in $R(\Psi)$).
\end{itemize}
The laws for identity have been taken care of by stipulation and associativity follows from Lemma \ref{depassociative}.
$R(\Phi)$ is clearly a finite category since $R(\Psi)$ is one, and it remains an inverse category since the only non-identity arrows that have been added to $R(\Phi)$ have codomain $K$. So we are done.
\end{construction}

\begin{problem}
To construct a function $S \colon FIC \rightarrow \FOLDS$.
\end{problem}
\begin{construction}\label{fictosig}
Let $\L$ be a FIC. We proceed by induction on $\oo{}{}$ and $\om{}{}$ to define a derivable signature expression consisting of a list of expressions $s(K_1),\dots,s(K_n)$ where $\oo{K_1}{K_2}\dots\oo{}{K_n} \in \ob{\L}$ and set $S(\L)=s(K_1),\dots,s(K_n)$. If $\ob{\L} = \varnothing$ then $\L = \varnothing$ and we let $S(\L)=\bullet$. If $\ob{\L} \neq \varnothing$ then take $K=K_m \in \ob{\L}$ and assume that $s(K_j)$ is defined for all $j<m$ such that $S = s(K_1),\dots,s(K_{m-1})$ is derivable. 
If $\cosieve{K}{\L} = \varnothing$ then we set $s(K) = K \colon \bullet$. 
Noting that $K \freshfor S$ since only $K_j$ for $j<m$ appear in $S$ by definition, we have the following derivation:
\[
\inferrule{
S \vdash
}
{
\inferrule
{
S; \bullet \vdash \\ K \freshfor S
}
{
S, K \colon \bullet \vdash
}
\:\: \ottdrulename{sig\_ext}
}
\:\: \ottdrulename{ctx\_empty}
\]
On the other hand, if $\cosieve{K}{\L} \neq \varnothing$ then let $K_1,\dots,K_m$ be be such that $\exists f \in \L(K,K_i)$ such that $f=gh \Rightarrow g=f$ or $h=f$. Since for each $i$ we have $\oo{K_i}{K}$ we can assume that $s(K_i)$ has been defined and that $S = s(K_1),\dots,s(K_m)$ is derivable, and we set $s(K) = T_K$. It now remains to check that $S, K \colon T_K$ is derivable. We proceed by $\om{}{}$-induction on the context $T_K$ which contains as variables exactly the elements of $\cosieve{K}{\L}$. So take some $f \colon K \rightarrow K'$ in $\cosieve{K}{\L}$ and assume that the context $\Gamma$ containing all $\om{g}{f}$ has been derived, which in particular means that $K' \colon T_{K'}$ is in $S$. Write $p_1,\dots,p_k$ for the arrows in $\cosieve{K'}{\L}$ and write $T_{K'}^-$ for the context $T_{K'}$ with the last declared variable removed, i.e. such that
\[
T_{K'} = T_{K'}^-,  p_k \colon K_{p_k} (s_{1} p_k,\dots,s_{l}p_k)
\]
where $s_{1},\dots,s_{l}$ are the (ordered) arrows in the cosieve $\cosieve{K_{p_k}}{\L}$.
Then we have
\[
\inferrule
{
\inferrule
{
S \vdash (p_1f,\dots,p_{k-1}f) \colon \Gamma \Rightarrow T_{K'}^- \\ S; \Gamma \vdash p_kf \colon K_{p_k} (s_{1} p_kf,\dots,s_{l}p_kf)
}
{
(K' \colon T_{K'}) \in S \\ S \vdash (p_1f,\dots,p_kf) \colon \Gamma \Rightarrow T_{K'}^-, p_k \colon K_{p_k} (s_{1} p_k,\dots,s_{l}p_k)
}
\:\: \ottdrulename{sub\_ext}
}
{
\inferrule
{
S; \Gamma \vdash K' (p_1f,\dots,p_kf) \colon \U
}
{
S; \Gamma, f \colon K' (p_1f,\dots,p_kf) \vdash
}
\:\: \ottdrulename{ctx\_ext}
}
\:\: \ottdrulename{type\_sort}
\]
But now note that by the inductive assumption, since $\om{p_kf}{f}$ we have that $\Gamma$ will contain $p_kf \colon K_{p_k} (s_{1} p_kf,\dots,s_{l}p_kf)$ and hence $S; \Gamma \vdash p_kf \colon K_{p_k} (s_{1} p_kf,\dots,s_{l}p_kf)$ is derivable (by \ottdrulename{type\_var}). 
On the other hand, $T_K^-$ is a context smaller than $T_K$ and therefore by induction on $k$ we can show that $S \vdash (p_1f,\dots,p_{k-1}f) \colon \Gamma \Rightarrow T_{K'}^-$, by repeating the argument given above.
\end{construction}

\begin{theorem}\label{sigbijfic}
$R$ and $S$ define a bijection $\FOLDS \simeq FIC$.
\end{theorem}
\begin{proof}
We proceed by induction to show that both $SR \colon \FOLDS \rightarrow \FOLDS$ and $SR \colon FIC \rightarrow FIC$ are equal to the identity.

Clearly, $SR(\bullet)\equiv\bullet$ and $RS(\varnothing)=\varnothing$.

Now assume that $SR(\Psi)\equiv\Psi$ and consider $\Phi \eqdef \Psi, K \colon \Gamma$.
Then by Construction \ref{fictosig} (and using the same notation as there) we have 
\[
SR(\Phi) \equiv s(R(\Psi)_1),\dots,s(R(\Psi)_m), K \colon T_K
\]
where $R(\Psi)_i$ are the objects of $R(\Psi)$, ordered in the usual way.
But by the inductive hypothesis we have 
$
\Psi \equiv S(R(\Psi))\equiv s(R(\Psi)_1),\dots,s(R(\Psi)_m)
$
and so it remains to check that $\Gamma \equiv T_K$. But note that by definition, $T_K$ is simply the total cosieve on $K$, and since $K$, as an object of $R(\Phi)$, is exactly the codomain of (the arrows represented by) the variables in $\Gamma$, we get that $T_K \equiv \Gamma$.

On the other hand, let $\L \in FIC$, $K_1,\dots,K_m$ its objects, ordered by $\oo{}{}$.
We proceed by induction on $m$ (i.e. on $\oo{}{}$), writing $\L'$ for the full subcategory of $\L$ consisting only of the objects $K_i$ with $i < m$.
By Construction \ref{fictosig} we have $S(\L) \equiv s(K_1), \dots, s(K_m)$ and therefore
\begin{align*}
RS(\L) &= R(s(K_1), \dots, s(K_m)) \\
	   &= R(s(K_1),\dots,s(K_{m-1})),K_m \colon T_{K_m} \\
\end{align*}
So on the assumption that $R(s(K_1),\dots,s(K_{m-1})) = \L'$ we have:
\begin{itemize}
\item $\ob{RS(\L)} = \ob{RS(\L')} \amalg \lbrace K \rbrace = \ob{\L}$
\item 
\[
RS(\L)(X,Y) = \fourpartdef{\morph{RS(\L')}(X,Y) = \L'(X,Y)}{ X,Y \neq K}{1_K}{X=Y=K}{\varnothing}{Y=K,X\neq K}{\vars{T_K}_Y = \L(K,Y) }{X=K,Y\neq K}
\]
Hence, $RS(\L)(X,Y) = \L(X,Y)$.
\item 
\[
g \catcomp_{RS(\L)} f = \threepartdef{g \catcomp_{RS(\L')}f = g \catcomp_{\L'} f}{X,Y, Z \neq K}{g}{f=1_K}{\dep{g}{f}}{X=K, Y,Z \neq K}
\]
\end{itemize}
This means exactly that $RS(\L)=\L$ and we are done.
\end{proof}

\bibliographystyle{shortalphabetic}

\bibliography{FICSasDTSrefs}

\begin{bibdiv}
\begin{biblist}

\bib{altenkirch16}{article}{
      author={Altenkirch, Thorsten},
      author={Kaposi, Ambrus},
       title={Type theory in type theory using quotient inductive types},
organization={ACM},
        date={2016},
     journal={ACM SIGPLAN Notices},
      volume={51},
      number={1},
       pages={18\ndash 29},
}

\bib{uniFOLDS}{article}{
      author={Ahrens, B.},
      author={North, P.},
      author={Shulman, M.},
       title={Univalent {F}{O}{L}{D}{S}},
        date={2014},
        note={Abstract available at:
  \url{http://www.ams.org/amsmtgs/2223_abstracts/1103-18-196.pdf}},
}

\bib{Cart86}{article}{
      author={Cartmell, J.},
       title={Generalized algebraic theories and contextual categories},
        date={1986},
     journal={Annals of Pure and Applied Logic},
      volume={32},
       pages={209\ndash 243},
  eprint={http://www.sciencedirect.com/science/article/pii/0168007286900539},
}

\bib{MFOLDS}{unpublished}{
      author={Makkai, M.},
       title={First order logic with dependent sorts with applications to
  category theory},
        date={1995},
  note={\url{http://www.math.mcgill.ca/makkai/folds/foldsinpdf/FOLDS.pdf}},
}

\bib{PalmgrenFOLDS}{article}{
      author={Palmgren, Erik},
       title={Categories with families, folds and logic enriched type theory},
        date={2016},
     journal={arXiv preprint arXiv:1605.01586},
}

\bib{PittsNS}{book}{
      author={Pitts, Andrew~M},
       title={Nominal sets: Names and symmetry in computer science},
   publisher={Cambridge University Press},
        date={2013},
}

\bib{shulman2015reedy}{article}{
      author={Shulman, Michael},
       title={Reedy categories and their generalizations},
        date={2015},
     journal={arXiv preprint arXiv:1507.01065},
}

\bib{TsemHMT}{article}{
      author={Tsementzis, D.},
       title={{Homotopy Model Theory I: Syntax and Semantics}},
        date={2016},
     journal={arXiv preprint arXiv:1603.03092},
}

\bib{TsemHSIP}{article}{
      author={Tsementzis, D.},
       title={{A Higher Structure Identity Principle}},
        date={2017},
     journal={arXiv preprint arXiv:1702.07776},
}

\bib{TsemWeavTTI}{article}{
      author={Weaver, M.},
      author={Tsementzis, D.},
       title={{Unfolding FOLDS}},
        note={In Preparation},
}

\end{biblist}
\end{bibdiv}

\end{document}